\documentclass[12pt,a4paper]{article}
\usepackage{amssymb}

\usepackage{graphicx,amssymb,amsfonts,epsfig,amsthm,a4,amsmath,url}
\usepackage[latin1]{inputenc}

\usepackage{enumitem}

\usepackage{amssymb}

\usepackage{graphicx,amssymb,amsfonts,epsfig,a4,amsmath,amsthm}
\usepackage[latin1]{inputenc}

\newtheorem{Thm}{Theorem}[section]
\newtheorem{Prop}[Thm]{Proposition}
\newtheorem{Lem}[Thm]{Lemma}
\newtheorem{Cor}[Thm]{Corollary}

\newtheorem{Rem}[Thm]{Remark}

\newcommand{\R}{\mathbf{R}}

\newcommand{\bprr}{\noindent \textbf{Proof} }
\newcommand{\epr}{~$\blacksquare$}

\title{Wasserstein distance and metric trees}
\author{Maxime MATHEY-PREVOT and Alain VALETTE}

\setlistdepth{9}
\newlist{myEnumerate}{enumerate}{9}
\setlist[myEnumerate,1]{label= }
\setlist[myEnumerate,2]{label= }
\setlist[myEnumerate,3]{label= }
\setlist[myEnumerate,4]{label= }
\setlist[myEnumerate,5]{label= }
\setlist[myEnumerate,6]{label= }
\setlist[myEnumerate,7]{label= }
\setlist[myEnumerate,8]{label= }
\setlist[myEnumerate,9]{label= }

\begin{document}

\maketitle

\hfill{\it A la m\'emoire de Vaughan, ami et bon vivant}
\begin{abstract} We study the Wasserstein (or earthmover) metric on the space $P(X)$ of probability measures on a metric space $X$. We show that, if a finite metric space $X$ embeds stochastically with distortion $D$ in a family of finite metric trees, then $P(X)$ embeds bi-Lipschitz into $\ell^1$ with distortion $D$. Next, we re-visit the closed formula for the Wasserstein metric on finite metric trees due to Evans-Matsen \cite{EvMat}. We advocate that the right framework for this formula is real trees, and we give two proofs of extensions of this formula: one making the link with Lipschitz-free spaces from Banach space theory, the other one algorithmic (after reduction to finite metric trees). 
\end{abstract}

%{\Large

\section{Introduction}

Embeddings of metric spaces, especially discrete metric spaces like graphs, into the Banach spaces $\ell^1$ or $L^1$, form a well-established part of metric geometry, with applications ranging from computer science to topology: we refer to \cite{NaorICM}, part I of \cite{DeLa} or Chapter 1 in \cite{Ostro}. In this paper we will be concerned with embeddings of Wasserstein spaces, that we now recall.

Let $(X,d)$ be a metric space and let $P_1(X)$ be the space of probability measures $\mu$ on $X$ with finite first moment, i.e. 
$$\int_X d(x_0,x)\,d\mu(x)<+\infty$$
for some (hence any) base-point $x_0\in X$. For compact $X$, the space $P_1(X)$ coincides with the space $P(X)$ of all probability measures on $X$.

The {\it Wasserstein metric} is a distance function on $P_1(X)$. Intuitively, given $\mu,\nu\in P_1(X)$, the distance $Wa(\mu,\nu)$ represents the amount of work necessary to transform $\mu$ into $\nu$. More precisely, a probability measure $\pi\in P(X\times X)$ is a {\it coupling} between $\mu$ and $\nu$ if its marginals are $\mu$ and $\nu$, i.e. $\mu(A)=\pi(A\times X)$ and $\nu(A)=\pi(X\times A)$ for any Borel subset $A\subset X$. And the Wasserstein distance $Wa(\mu,\nu)$ is defined as
$$Wa(\mu,\nu)=\inf\Big\{\int_{X\times X}d(x,y)\,d\pi(x,y):{\mbox{$\pi$ coupling between $\mu$ and $\nu$}}\Big\}.$$
Note that $X$ embeds isometrically in $P_1(X)$ by $x\mapsto \delta_x$ (the Dirac mass at $x$). See Chapter 5 in \cite{Santa} or Chapter 7 in \cite{Villa} for more on the Wasserstein distance, also called {\it Kantorovitch-Rubinstein distance} or {\it earthmover distance} (EMD) in computer science papers. We denote by $Wa(X)$ the space $P_1(X)$ endowed with the Wasserstein distance, and call it {\it the Wasserstein space} of $X$. For a coupling $\pi$, the \textit{cost} of $\pi$ is the quantity $\int_{X \times X} d(x,y) d\pi(x,y)$.

Let $\mathcal{Y}=(Y_i,d_i)_{i\in I}$ be a finite family of metric spaces. We say that a metric space $(X,d)$ {\it embeds stochastically in $\mathcal{Y}$ with distortion $D\geq 1$} if there exists non-negative numbers $(p_i)_{i\in I}$ summing up to 1, and maps $f_i:X\rightarrow Y_i$ (for each $i\in I$) such that:
\begin{itemize}
\item Each $f_i$ is non-contracting, i.e. for every $x,y\in X$ we have $d_i(f_i(x),f_i(y))\geq d(x,y)$.
\item For every $x,y\in X$ we have $\sum_{i\in I} p_i d_i(f_i(x),f_i(y))\leq D\cdot d(x,y).$
\end{itemize}

The first aim of this paper is to prove the following result; 

\begin{Thm}\label{stoch} Assume that the finite metric space $(X,d)$ embeds stochastically with distortion $D$ into a family of finite metric trees. Then $Wa(X)$ embeds bi-Lipschitz into $\ell^1$ with distortion at most $D$.
\end{Thm}

Here, by a metric tree, we mean a tree $T=(V,E)$ endowed with a positive weight function $w:E\rightarrow\R_{>0}:e\mapsto w_e$. For $x,y\in V$ we denote by $[x,y]$ the set of edges on the unique path from $x$ to $y$ and we endow $V$ with the distance $d_T(x,y)=\sum_{e\in [x,y]} w_e$.

We learned Theorem \ref{stoch} from the paper \cite{InTha} by P. Indyk and N. Thaper, who get a less precise $O(D)$ for the distortion of the embedding into $\ell^1$, and provide a rather frustrating comment that prompted our desire to provide a direct proof of Theorem \ref{stoch}.\footnote{We provide the comment for completeness: ``{\it The embedding can be seen as resulting from a combination of the following two results:

1. The result of \cite{Cha}, who (implicitely) showed that the techniques of \cite{KleTar} imply the following: if a metric $M$ can be probabilistically embedded into trees with distortion $c$, then the EMD over $M$ can be embedded into $\ell^1$ with distortion $O(c)$.

2. The result of \cite{CCGGP} who showed that the Euclidean metric over $\{1,...,\Delta\}^d$ can be probabilistically embedded into trees with distortion $O(d\log \Delta)$. Again that result is implicit in that paper.}''}

It was shown by J. Fakcharoenphol, S. Rao and K. Talwar \cite{FaRaTa} that any finite metric space on $n$ points embeds stochastically with distortion $O(\log n)$ into a family of finite metric trees (and this bound is optimal). Using this it was shown by F. Baudier, P. Motakis, G. Schlumprecht and A. Z\`sak (Corollary 8 in \cite{BMSZ}) that, for $X$ a finite metric space on $n$ points, the lamplighter metric space $La(X)$ embeds into $\ell^1$ with distortion $O(\log n)=O(\log\log|La(X)|)$. Using the same result from \cite{FaRaTa}, our Theorem \ref{stoch} immediately implies:

\begin{Cor} For any finite metric space $X$ on $n$ points, the Wasserstein space $Wa(X)$ embeds bi-Lipschitz into $\ell^1$ with distortion at most $O(\log n)$.
\epr
\end{Cor}

Combining with the isometric embedding $X\rightarrow Wa(X):x\mapsto \delta_x$, we get as corollary a celebrated result by J. Bourgain \cite{Bour} \footnote{Of course all the difficulty becomes hidden in \cite{FaRaTa}!}

\begin{Cor} Any finite metric space on $n$ points, embeds bi-Lipschitz into $\ell^1$ with distortion at most $O(\log n)$.
\end{Cor}

It turns out that on finite metric trees there is a remarkable closed formula for the Wasserstein distance. It originated in papers in computer science in 2002 and probably earlier: see Charikar \cite{Cha}, for measures supported on the leaves of the tree\footnote{A proof for this special case appears in lemma 3.1 of \cite{Kloeck}.}. For general probability measures on a finite metric tree, the formula appears in a paper in biomathematics (see section 2 in S.N. Evans and F.A. Matsen \cite{EvMat}). We believe it deserves to be better known in mathematical circles. To understand it, let $T=(V,E)$ be a metric tree, fix a base-vertex $x_0\in V$ (so that $T$ appears as a rooted tree). Any edge $e\in E$ separates $T$ into two half-trees, and we denote by $T_e$ the set of vertices of the half-tree NOT containing $x_0$: if we view the tree as hanging from the root, $T_e$ is the subtree hanging below the edge $e$.  

\begin{Thm}\label{WassTrees} Let $T=(V,E)$ be a finite, rooted metric tree. Then for $\mu,\nu\in P(V)$:
\begin{equation}\label{LAforrmule!}
Wa(\mu,\nu)=\sum_{e\in E} w_e|\mu(T_e)-\nu(T_e)|
\end{equation}
\end{Thm}

This formula has numerous implications: first, the RHS is independent of the choice of the root; second, it shows that the Wasserstein metric on $T$ is a $L^1$-metric (see lemma \ref{treesinL1} below). 

Our second aim in this paper is to give two new proofs of Theorem \ref{WassTrees}. The first one advocates that the right framework for Theorem \ref{WassTrees} is real trees: by exploiting a connection with the theory of Lipschitz-free spaces from Banach space theory, we will extend the result to metric trees with countably many vertices. The second proof is by double inequality: the inequality $Wa(\mu,\nu)\geq\sum_{e\in E} w_e|\mu(T_e)-\nu(T_e)|$ follows by considering the canonical embedding of the tree into $\ell^1$ and its barycentric extension to $P_1(V)$. The converse inequality is proved by first reducing to finite metric trees and, for those, given $\mu,\nu\in P(V)$, by providing an algorithmic construction of a coupling $\pi$ with $\int_{V\times V}d(x,y)\,d\pi(x,y)=\sum_{e\in E} w_e|\mu(T_e)-\nu(T_e)|$. 

%\medskip
%DIRE POURQUOI ON S'INTERESSE AUX PLONGEMENTS DANS $L^1$!

\medskip
The paper is organized as follows. In section 2 we prove Theorem \ref{stoch}, taking Theorem \ref{WassTrees} for granted. Sections 3 and 4 present our two proofs of Theorem \ref{WassTrees}, suitably generalized to metric trees with countably many vertices (see Theorem \ref{Wasscountabletree}). Finally the Appendix provides a comparison between various $\sigma$-algebras of sets on a real tree, that appeared in the literature.

\medskip
{\bf Acknowledgements}: We thank Florent Baudier for many insightful discussions, for providing us with lemma \ref{embeddingWass} and its proof, and for useful comments on an earlier version of the paper; and Gilles Lancien for pointing out the connection with Lipschitz-free spaces.

\section{Stochastic embeddings}

We will prove Theorem \ref{stoch} by means of a series of lemmas. The first one is lemma 3 in \cite{BMSZ} .

\begin{Lem}\label{composing} Assume that the metric space $(X,d)$ embeds stochastically in $\mathcal{Y}=(Y_i,d_i)_{i\in I}$, and that for every $i\in I$ the space $Y_i$ embeds bi-Lipschitz into $\ell^1$ with distortion $C_i\geq 1$. Then $X$ embeds bi-Lipschitz into $\ell^1$ with distortion at most $CD$, where $C=\max_{i\in I} C_i$. \hfill\epr
\end{Lem}

The second lemma was suggested to us by F. Baudier.

\begin{Lem}\label{embeddingWass} If the finite metric space $(X,d)$ embeds stochastically into $\mathcal{Y}=(Y_i,d_i)_{i\in I}$ with distortion $D$, then $Wa(X)$ embeds stochastically into $(Wa(Y_i))_{i\in I}$ with distortion $D$.
\end{Lem} 

\bprr For $i\in I$, let $p_i\geq 0$ and $f_i:X\rightarrow Y_i$ be realizing the stochastic embedding with distortion $D$ of $X$ into $\mathcal{Y}$. Consider then $(f_i)_*:P(X)\rightarrow P(Y_i):\mu\mapsto (f_i)_*(\mu)$. We claim that the stochastic embedding with distortion $D$ of $Wa(X)$ into the family $(Wa(Y_i))_{i\in I}$ is realized by the $p_i$'s and the $(f_i)_*$'s; to see this, we check the two points in the definition of a stochastic embedding. Fix $\mu,\nu\in Wa(X)$.
\begin{itemize}
\item Fix $i\in I$. Let $\pi_i$ be a coupling between $(f_i)_*(\mu)$ and $(f_i)_*(\mu)$ such that $Wa((f_i)_*(\mu),(f_i)_*(\nu))=\sum_{y,y'\in Y_i}d_{Y_i}(y,y')\pi_i(y,y')$. For $y\in Y_i\setminus f_i(X)$, we have 
$$\sum_{y'\in Y_i}\pi_i(y,y')=(f_i)_*(\mu)(y)=\mu(f_i^{-1}(y))=0$$
hence $\pi_i(y,y')=0$ for every $y'\in Y_i$. So $\pi_i$ vanishes outside of $f_i(X)\times f_i(X)$. Hence we may define $\sigma_i\in P(X\times X)$ by $\sigma_i(x,x')=\pi_i(f_i(x),f_i(x'))$ and $\sigma_i$ is a coupling between $\mu$ and $\nu$. Then 
$$Wa((f_i)_*(\mu),(f_i)_*(\nu))=\sum_{y,y'\in Y_i}d_{Y_i}(y,y')\pi_i(y,y')=\sum_{x,x'\in X} d_{Y_i}(f_i(x),f_i(x'))\pi_i(f_i(x),f_i(x'))$$
$$=\sum_{x,x'\in X} d_{Y_i}(f_i(x),f_i(x'))\sigma_i(x,x')\geq \sum_{x,x'\in X} d(x,x')\sigma(x,x')\geq Wa(\mu,\nu),$$
where the first inequality follows from the fact that $f_i$ is non-contracting. So $(f_i)_*$ is non-contracting as well.

\item Let $\pi\in P(X\times X)$ be a coupling between $\mu$ and $\nu$ such that $Wa(\mu,\nu)=\sum_{x,x'\in X} d(x,y)\pi(x,x')$. Set $\tau_i=(f_i\times f_i)_*(\pi) \in P(Y_i\times Y_i)$. Then $\tau_i$ is a coupling between $(f_i)_*(\mu)$ and $(f_i)_*(\nu)$ and
$$\sum_{i\in I} p_iWa((f_i)_*(\mu),(f_i)_*(\nu))\leq \sum_{i\in I} p_i\sum_{y,y'\in Y_i}d_{Y_i}(y,y)\tau_i(y,y')$$
$$=\sum_{i\in I} p_i\sum_{x,x'\in X} d_{Y_i}(f_i(x),f_i(x'))\tau_i(f_i(x),f_i(x'))=\sum_{i\in I} p_i\sum_{x,x'\in X} d_{Y_i}(f_i(x),f_i(x'))\pi(x,x')$$
$$=\sum_{x,x'\in X}\pi(x,x')\sum_{i\in I}p_id_{Y_i}(f_i(x),f_i(x'))\leq D\cdot\sum_{x,x'\in X}\pi(x,x')d(x,x')=D\cdot Wa(\mu,\nu),$$
where the second inequality follows from the fact that the $f_i$'s provide a stochastic embedding. This concludes the proof.
\end{itemize}
\epr

\medskip
Combining lemmas \ref{composing} and \ref{embeddingWass} we immediately get:

\begin{Cor}\label{immediate} If the finite metric space $(X,d)$ embeds stochastically into $\mathcal{Y}=(Y_i,d_i)_{i\in I}$ with distortion $D$, and each $Wa(Y_i)$ embeds bi-Lipschitz in $\ell^1$ with distortion $C_i$, then $Wa(X)$ embeds bi-Lipschitz into $\ell^1$ with distortion at most $CD$, where $C=\max_{i\in I} C_i$.
\hfill\epr
\end{Cor}

To prove Theorem \ref{stoch}, in view of Corollary \ref{immediate}, it is therefore enough to observe:

\begin{Lem}\label{treesinL1} If $T=(V,E)$ is any finite metric tree, then $Wa(T)$ embeds isometrically into $\ell^1$.
\end{Lem} 

\bprr Fix a root $x_0\in V$ and, for any edge $e\in E$, let $T_e\subset V$ be defined as in the Introduction. The map 
$$Wa(T)\rightarrow \ell^1(E):\mu\mapsto (e\mapsto w_e\mu(T_e))$$
is an isometric embedding of $Wa(T)$, by Theorem \ref{WassTrees}.
\epr

\medskip
This concludes the proof of Theorem \ref{stoch} (taking Theorem \ref{WassTrees} for granted).

\section{First proof of Theorem \ref{WassTrees}}

\subsection{Lipschitz-free spaces}

For a metric space $(X,d)$ with a base-point $x_0\in X$, we denote by $Lip_0(X)$ the Banach space of Lipschitz functions on $X$ vanishing at $x_0$, endowed with the Lipschitz norm. 
The space $Lip_0(X)$ has a canonical pre-dual, called the {\it Lipschitz-free space} of $X$ (see e.g. Chapter 2 in \cite{Weaver}, Chapter 10 in \cite{Ostro}) and denoted by $\mathcal{F}(X)$: it is the closed linear subspace of the dual space $Lip_0(X)^*$ generated by the point evaluations $\delta_x\;(x\in X\setminus\{x_0\})$.

For $\mu\in P_1(X)$, the linear form $f\mapsto \int_X f(x)\,d\mu(x)$ defines an element of the dual $Lip_0(X)^*$: this way we get an embedding of $Wa(X)$ into $Lip_0(X)^*$. When $X$ is a complete separable metric space, it can be shown that this is actually  an isometric embedding of $Wa(X)$ into $\mathcal{F}(X)$ (see Theorem 1.13 in \cite{OsOs} or section 2 in \cite{NaSche}).%\footnote{The authors of \cite{OsOs} write: {\it ``We believe that this result is known to experts. However, since a suitable reference has not been found, its proof is presented.''}})

\subsection{Real trees}

A {\it real tree} $(T,d)$ is a geodesic metric space which is 0-hyperbolic in the sense of Gromov. For $x,y\in T$, we denote by $[x,y]$ the {\it segment} between $x$ and $y$, i.e. the unique arc joining them. A point $x\in T$ is a {\it branching point} if $T\setminus\{x\}$ has at least 3 connected components; we denote by $Branch(T)$ the set of branching points of $T$. Fix a base-point $x_0\in T$. For $x\in T$, we set
$$T_x=\{y\in T: x\in [x_0,y]\}=\{y\in T: d(x_0,y)=d(x_0,x)+d(x,y)\};$$
so letting $T$ hang from the root $x_0$, the set $T_x$ is the part of $T$ lying below $x$. 

Following A. Godard \cite{God}) we say that a subset $A\subset T$ is {\it measurable} if, for every $x,y\in T$, the set $A\cap [x,y]$ is Lebesgue-measurable in $[x,y]$. On the $\sigma$-algebra $\mathcal{G}$ of measurable subsets, there is a unique measure $\lambda$ such that $\lambda([x,y])=d(x,y)$: we call $\lambda$ the {\it length measure} \footnote{ See the Appendix below for a comparison of various $\sigma$-algebras associated with real trees.}. It is defined as follows: for $S$ a segment in $T$, let $\lambda_S$ denote Lebesgue measure on $S$. Then, for $A\in \mathcal{G}$, if $R\subset T$ is a finite disjoint union of segments, say $R=\cup_{i=1}^k S_i$, we set $\lambda_R(A)=\sum_{i=1}^k \lambda_{S_i}(S_i\cap A)$. Finally we set 
\begin{equation}\label{Godmeasure}
\lambda(A)=\sup_{R\in\mathcal{R}} \lambda_R(A)
\end{equation}
where $\mathcal{R}$ is the set of subsets of $T$ that can be expressed as finite disjoint unions of segments.

For $A$ a closed subset of $T$ containing $x_0$, still following \cite{God} we define a function $L_A: A\rightarrow\R^+$ by $L_A(a)=\inf\{d(a,x):x\in A\cap [x_0,a[\}$. So $L_A(a)>0$ if and only if $a$ is isolated in $A\cap [x_0, a]$. We then define a measure $\mu_A$ on $A$ by $\mu_A=\lambda|_A + \sum_{a\in A} L_A(a)\delta_a$. In Theorem 3.2 in \cite{God}, it is proved that, if $A$ is a closed subset of $T$ containing $Branch(T)$, then $\mathcal{F}(A)$ is isometrically isometric to $L^1(A,\mu_A)$. 

Assume from now on that the real tree $T$ is complete and separable. Then by the previous sub-section, for $A$ a closed subset of $T$ containing $Branch(T)$, the space $Wa(A)$ isometrically embeds into $L^1(A,\mu_A)$. This embedding is not written explicitly in \cite{God}; by making it explicit we get a closed formula for the Wasserstein distance on closed subsets of real trees.

\begin{Prop}\label{Wassrealtrees} Let $(T,d)$ be a complete, separable real tree and let $A$ be a closed subset of $T$ containing $Branch(T)$. For $\mu,\nu\in Wa(A)$ we have:
\begin{equation}\label{eq1}
Wa(\mu,\nu)=\int_A |\mu(T_x\cap A)-\nu(T_x\cap A)|\,d\mu_A(x).
\end{equation}
\end{Prop}

\bprr By the proof of Theorem 3.2 in \cite{God}, the map
$$\Phi: L^\infty(A,\mu_A)\rightarrow Lip_0(A):g\mapsto \Big(a\mapsto \int_{[x_0,a]\cap A} g(x)\,d\mu_A(x)\Big)$$
is an isometric isomorphism which is weak$^*$-weak$^*$ continuous, so its transpose $\Phi^*$ realizes the desired isometric isomorphism $\mathcal{F}(A)\rightarrow L^1(A,\mu_A)$. Denoting by $\chi_{[x,y]}$ the characteristic function of the interval $[x,y]$, the previous formula may be re-written:
$$(\Phi(g))(a)=\int_A \chi_{[x_0,a]}(x)g(x)\,d\mu_A(x).$$
For $\nu\in Wa(A)$, we compute $\Phi^*(\nu)$. For $g\in L^\infty(A,\mu_A)$, we have
$$(\Phi^*(\nu), g)=(\nu,\Phi(g))=\int_A (\Phi(g))(a)\,d\nu(a)=\int_A\Big(\int_A \chi_{[x_0,a]}(x)g(x)\,d\mu_A(x)\Big)\,d\nu(a).$$
As the measure $\mu_A$ is $\sigma$-finite (here we use that the tree $T$ is separable), we may appeal to Fubini:
$$(\Phi^*(\nu), g)=\int_A g(x)\Big(\int_A\chi_{[x_0,a]}(x)\,d\nu(a)\Big)\,d\mu_A(x)=\int_A g(x)\nu(T_x\cap A)\,d\mu_A(x).$$
Since this holds for every $g\in L^\infty(A,\mu_A)$ we deduce that, for almost every $x\in A$:
$$(\Phi^*(\nu))(x)=\nu(T_x\cap A).$$
Equation (\ref{eq1}) follows.
\epr

\medskip
\begin{Rem} When $A=T$, Proposition \ref{Wassrealtrees} becomes, for $T$ a complete separable real tree and $\mu,\nu\in Wa(T)$:
\begin{equation}\label{eq2}
Wa(\mu,\nu)=\int_T |\mu(T_x)-\nu(T_x)|\,d\lambda(x).
\end{equation}
When $T$ is the geometric realization of a finite metric tree, equation (\ref{eq2}) appears as equation (5) in \cite{EvMat}; the proof is different.
\end{Rem}

\begin{Thm}\label{Wasscountabletree} Let $T=(V,E)$ be a rooted metric tree with countably many vertices. Then for $\mu,\nu\in P_1(V)$:
$$Wa(\mu,\nu)=\sum_{e\in E} w_e|\mu(T_e)-\nu(T_e)|$$
\end{Thm}

\bprr Fix $\mu,\nu\in P_1(V)$. For an edge $e$, let $e^+,e^-$ be the vertices of $e$, chosen so that $d(x_0,e^+) < d(x_0,e^-)$. Then on the arc $[e^+,e^-]$ the function $x\mapsto |\mu(T_x)-\nu(T_x)|$ is constant, equal to $|\mu(T_e)-\nu(T_e)|$. So by formula (\ref{eq2}):
$$Wa(\mu,\nu)=\sum_{e\in E}\int_{[e^+,e^-]} |\mu(T_x)-\nu(T_x)|\,d\lambda(x)= \sum_{e\in E} w_e|\mu(T_e)-\nu(T_e)|.$$
\epr

%Let $T=(V,E,\ell)$ be a locally finite weighted tree, meaning that the graph $(V,E)$ is a locally finite tree, and $w:E\rightarrow\R_{>0}:e\mapsto w_e$ is a positive weight function, that we view as a positive measure on $E$. For $x,y\in V$ we denote by $[x,y]$ the set of edges on the unique path from $x$ to $y$. We endow $V$ with the distance
%$$d(x,y)=\sum_{e\in [x,y]} w_e,$$
%and we assume that $(V,d)$ is a proper metric space, i.e. balls are finite. This is the case if $w$ is bounded below by a positive constant, e.g. for the combinatorial distance on $T$ (i.e. $w \equiv 1$).

%Fix a base-vertex $x_0\in V$,  Any edge $e\in E$ divides $T$ in two subtrees, by removing it. We denote by $T_e$ the vertex set of the subtree NOT containing $x_0$.

\section{Second proof of Theorem \ref{WassTrees}}

\subsection{Lipschitz maps to Banach spaces}

%Let $(X,d)$ be a metric space. We denote by $P_1(X)$ the space of probability measures $\mu$ on $X$ such that 
%$$\int_X d(x_0,x)\,d\mu(x)<+\infty$$
%for some (hence any) $x_0\in X$. As in \cite{Villa}, we endow $P_1(X)$ with the {\bf Wasserstein distance}:
%$$Wass(\mu,\nu)=\inf_{\pi\in\Pi(\mu,\nu)} \int_{X\times X} d(x,y)\,d\pi(x,y),$$
%where $\Pi(\mu,\nu)$ denotes the set of probability measures $\pi$ on $X\times X$ having $\mu$ and $\nu$ as marginal distributions: 

\begin{Prop}\label{extension} Let $(X,d)$ be a metric space, and let $E$ be a Banach space. Any $C$-Lipschitz map $\beta:X\rightarrow E$ extends canonically to a $C$-Lipschitz map $\tilde{\beta}:Wa(X)\rightarrow E:\mu\mapsto\tilde{\beta}(\mu)$ defined as the barycenter of $\beta(X)$ with respect to $\mu$, i.e.
$$\tilde{\beta}(\mu)=\int_X \beta(x)\,d\mu(x).$$
\end{Prop}

\bprr Let $x_0$ be a base-point in $X$. Composing $\beta$ with a translation in $E$, we may assume that $\beta(x_0)=0$. Then, as $\|\beta(x)-\beta(y)\|\leq C\cdot d(x,y)$ for any $x,y\in X$, we get $\|\beta(x)\|\leq C\cdot d(x_0,x)$, hence $\|\tilde{\beta}(\mu)\|\leq\int_X\|\beta(x)\|\,d\mu(x)\leq C\int_X d(x_0,x)\,d\mu(x)<+\infty$. So $\tilde{\beta}$ is well-defined.

To check that $\tilde{\beta}$ is $C$-Lipschitz, observe that for $\mu,\nu\in P_1(X)$ and $\pi$ a coupling between $\mu$ and $\nu$, we have:
$$\|\tilde{\beta}(\mu)-\tilde{\beta}(\nu)\|=\Big\|\int_X\beta(x)\,d\mu(x)-\int_X\beta(y)\,d\nu(y)\Big\|=\Big\|\int_{X\times X}\beta(x)\,d\pi(x,y)-\int_{X\times X}\beta(y)\,d\pi(x,y)\Big\|$$
$$\leq \int_{X\times X}\|\beta(x)-\beta(y)\|\,d\pi(x,y)\leq C\int_{X\times X}d(x,y)\,d\pi(x,y).$$
The result follows by taking the infimum over all couplings $\pi$.
\epr

\medskip
\begin{Rem}
Observe that, if $\beta$ in Proposition \ref{extension} is bi-Lipschitz, in general its extension $\tilde{\beta}$ is not. Indeed take $E=\R$, and let $X\subset\R$ be any subset with at least 3 elements, the inclusion $\beta:X\rightarrow\R$ is isometric, but $\tilde{\beta}$ is not even injective. 
\end{Rem}

\medskip
Let $T=(V,E)$ be a metric tree; we denote by $\chi_{[x,y]}$ the characteristic function of the set of edges in $[x,y]$. There is a well-known isometric embedding $\beta: V\rightarrow\ell^1(E,w):x\mapsto\chi_{[x_0,x]}$
(it is hard to locate the first appearance of this embedding in the literature: we learned it from \cite{Haagerup}). By Proposition \ref{extension}, we extend it to a 1-Lipschitz map $\tilde{\beta}:P_1(V)\rightarrow\ell^1(E,w)$. Ultimately we will see that $\tilde{\beta}$ is isometric. For the moment we prove:

\begin{Prop}\label{trees} Let $T=(V,E)$ be a metric tree. For $\mu,\nu\in P_1(V)$:
$$\|\tilde{\beta}(\mu)-\tilde{\beta}(\nu)\|_1=\sum_{e\in E}w_e|\mu(T_e)-\nu(T_e)|\leq Wa(\mu,\nu).$$
\end{Prop}

\bprr The inequality follows from Proposition \ref{extension}, we focus on the equality. But 
$$\|\tilde{\beta}(\mu)-\tilde{\beta}(\nu)\|_1=\sum_{e\in E}w_e |\tilde{\beta}(\mu)(e)-\tilde{\beta}(\nu)(e)|.$$
So it it enough to prove that $\tilde{\beta}(\mu)(e)=\mu(T_e)$. So we compute:
$$\tilde{\beta}(\mu)(e)=\sum_{x\in V}\beta(x)(e)\mu(x)=\sum_{x\in V}\chi_{[x_0,x]}(e)\mu(x)=\sum_{x\in T_e}\mu(x)=\mu(T_e)$$
as $\chi_{[x_0,x]}(e)=1$ if and only if $x\in T_e$.
\epr

\bigskip
Our aim now is to prove that the inequality in Proposition \ref{trees} is actually an equality, i.e for metric trees $T=(V,E)$ with countably many vertices we wish to prove the reverse inequality
\begin{equation}\label{ineq}
\sum_{e\in E}w_e|\mu(T_e)-\nu(T_e)|\geq Wa(\mu,\nu).
\end{equation}
%Since the set of finitely supported probability measures is dense in $Wa(V)$ (by Theorem 1.13 of \cite{OsOs}), to prove inequality (\ref{ineq}) we may reduce to finitely supported measures. In other words, we may restrict to finite metric trees, as the convex hull of a finite set is a finite metric tree.
%$(P_1(V), Wa)$ is metric and thus separated which implies that $P_1(V) \times P_1(V)$ is separated too (with the product topology).
Theorem 1.13 of \cite{OsOs} implies that the set of finitely supported probability measures is dense in $(P_1(V), Wa)$. Of course $Wa(\cdot , \cdot) : P_1(V) \times P_1(V) \to \R$ is continuous, and
$$P_1(V) \times P_1(V) \to \R: (\mu, \nu) \mapsto \sum_{e \in E} w_e | \mu(T_e) - \nu(T_e) |$$
 is continuous too, as an immediate consequence of Proposition \ref{trees}. So to show (\ref{ineq}) we may restrict for finitely supported measures i.e. we may restrict to finite metric trees. %To see that $(\mu, \nu) \mapsto \sum_{e \in E} w_e | \mu(T_e) - \nu(T_e) |$ is continuous we appeal to Theorem 7.12 in \cite{Villa} which states (among other things) that if $\mu \in P_1(V)$ and if $(\mu_k)$ is a sequence in $ P_1(V)$ then $Wa(\mu_k,\mu) \underset{k \to+\infty}{\longrightarrow} 0$ implies that $\int_V \phi d\mu_k \underset{k \to+\infty}{\longrightarrow} \int_V \phi d\mu$ for all $\phi$ on $V$ satisfying $|\phi(x) | \leq C(1+d(x_0,x))$ for some $C \in \R$. Precisely we take $(\mu_k, \nu_k)$ a sequence in $P_1(V) \times P_1(V)$ converging to $(\mu,\nu)$ in the product topology, then: %$Wa(\mu_k,\mu)  \underset{k \to+\infty}{\longrightarrow} 0$, $Wa(\nu_k,\nu)  \underset{k \to+\infty}{\longrightarrow} 0$, and:

\subsection{An algorithm for finite metric trees}

%If $T = (V,E)$ is a rooted-tree with root $x_0$ and if $e \in E$, we already defined $ T_e$ as the half-tree which set of vertices do not contain $x_0$. We have to use some more definitions: If $v,w \in V$ we write $[v,w]$ the set of vertices in the (unique) path from $v$ to $w$ in the tree. If $e \in E$ is an edge, we write $e^+$ and $e^-$ its two extremities, $e^+$ being closer to the root (i.e. $d(x_0,e^+) < d(x_0,e^-)$). If $v, w \in V$, we say that $w$ is a descendent of $v$ if  $v \in [x_0,w]$ (notice a vertex is its own descendent) and we say that $w$ is a child of $v$ - and that $v$ is the parent of $w$ - if $w$ is a descendent of $v$ and $[v,w] = \{w,v\}$. If $v \in V$ we write $T_v$ the half tree which set of vertices is the set of all descendent of $v$, hence $T_{x_0} = T$ and if $e \in E$ then $T_e = T_{e^-}$.\\

Let $T=(V,E)$ be a finite metric tree. Recall from the proof of theorem \ref{Wasscountabletree} that if $e \in E$ is an edge, we write $e^+$ and $e^-$ its two extremities chosen so $d(x_0,e^+) < d(x_0,e^-)$, moreover if $v, w \in V$, we say that $w$ is a descendant of $v$ if  $v \in [x_0,w]$ (notice that a vertex is its own descendant) and we say that $w$ is a child of $v$ - and that $v$ is the parent of $w$ - if $w$ is a descendent of $v$ and $[v,w] = \{w,v\}$. If $v \in V$ we write $T_v$ the half tree with set of vertices the set of all descendants of $v$, hence $T_{x_0} = T$ and if $e \in E$ then $T_e = T_{e^-}$.\\

To show that 
\begin{align*} 
Wa(\mu, \nu) \leq \sum_{e \in E} w_e \mid \mu(T_e) - \nu(T_e) \mid
\end{align*}
we provide an algorithm which transforms a probability measure $\mu'$, initially set to $\mu$ into $\nu$. In parallel, this algorithm keeps track of a variable (here a matrix) $\pi' = (\pi' (x,y))_{x,y \in V} := (\pi'_{x,y})_{x,y \in V}  $ that, all the way through the running of the algorithm, provides a coupling between $\mu$ and $\mu'$. When the algorithm stops we will have $\mu' = \nu$ and the cost of the coupling $\pi'$ will be $\sum_{e \in E} w_e \mid \mu(T_e) - \nu(T_e) \mid$. This algorithm runs in two \textit{phases}; intuitively speaking the first phase brings up (towards the root) the excess of mass from those subtrees $T_e$  with $\mu(T_e) > \nu(T_e)$, and the second phase let that mass fall (towards the leaves) in the subtrees $T_e$ with  $\mu(T_e) < \nu(T_e)$. Still intuitively, for every vertex $x$, $\pi'_{x,y}$ is the mass attributed by $\mu'$ to $x$ coming from $y$; the coupling remembers where the mass comes from. We consider that the vertices of $T$ are numbered with $1,2,...,n:=|V|$ (e.g. in such a way that given two vertices that are at distinct depths in the tree, the deeper one is associate to a lower number than the other). The algorithm is such that it moves first the mass coming from vertices with a low number.

\paragraph{Algorithm. } \textit{ }\\ 
\% \textit{\underline {Initialization}}:\\
$\mu ' \leftarrow \mu$.\\
\textbf{for all} $v$
\begin{myEnumerate}
\item $\pi_v' \leftarrow \vec{0}$
\item $\pi_{v,v}' \leftarrow \mu(v)$
\end{myEnumerate}
\textbf{end for}\\
$M \leftarrow 0$ \% \textit{This variable is used just for the proof}\\
\% \textit{\underline{Phase (1)}}:\\
\textbf{for} $N$ \textbf{depth level, from the deeper up to} 1:
\begin{myEnumerate}
\item \textbf{for all} $T_e$ \textbf{subtree whose root} $e^-$ \textbf{is at depth} $N$: \% \textit{Loop (*)}
\begin{myEnumerate}
%\item \textbf{Let} $e^+$ \textbf{be the parent of} $e^-$
\item \textbf{if} $\mu'(T_e) > \nu(T_e)$ \textbf{then}
\begin{myEnumerate}
\item \% \textit"{we bring $ (\mu'(T_e) - \nu(T_e)) $ up one level":} \\
$ x \leftarrow (\mu'(T_e) - \nu(T_e)) $  \\
$\mu'(e^-) \leftarrow \mu'(e^-) - x $  \% \textit{and simultaneously}  \\
 $\mu'(e^+)\leftarrow \mu'(e^+) + x $. \\
%Le coût de cette modification est $\ell_T (\mu'(T) - \nu(T)) =  \ell_T |\mu'(T) - \nu(T)|$.\\
$j \leftarrow  \min\{k : \sum_{i=1}^k \pi_{e^-,i}' \geq x \}$  \\
\textbf{for }$i < j$
\begin{myEnumerate}
\item $\pi_{e^+,i}' \leftarrow \pi_{e^+,i}' + \pi_{e^-,i}'$
\end{myEnumerate}
\textbf{end for}\\
$\pi_{e^+,j}'\leftarrow \pi_{e^+,j}' + (x - \sum_{n=1}^{j-1}\pi_{e^-,n})$\\
$\pi_{e^-,j}'\leftarrow \pi_{e^+,j}' - (x - \sum_{n=1}^{j-1}\pi_{e^-,n})$\\
\textbf{for all} $i < j$ 
\begin{myEnumerate}
\item $\pi_{e^-,i} \leftarrow 0$
\end{myEnumerate}
\textbf{end for}

\end{myEnumerate}
\textbf{end if}
\end{myEnumerate}
$M \leftarrow M+1$\\
\textbf{end for}
\end{myEnumerate}
\textbf{end for}\\
%%%%%%%%%%%%%%%%%%%%%%%%%%%%%%%%%%%%%%%%%%%%%%%%%%%%%%%%
\% \textit{\underline{Phase (2)}}:\\
\textbf{for} $N$ \textbf{depth level from }$0$\textbf{ to the deepest level in the tree}$ - 1$:
\begin{myEnumerate}
\item \textbf{for all} $T$ \textbf{ subtree whose root} $r$ \textbf{is at depth} $N$:
\begin{myEnumerate}
\item \textbf{let} $s_1,...,s_n$ \textbf{be the sons of }$r$ 
\item \textbf{if} $\mu'(r)> \nu(r)$:
\begin{myEnumerate} 
\item \textbf{for} $i=1... n$ \% \textit{Loop (**)}:
\begin{myEnumerate}
\item \textbf{if} $\nu(T_{s_i}) > \mu'(T_{s_i}) $:
\begin{myEnumerate}
\item \% \textit{"We let $(\nu(T_{s_i}) - \mu'(T_{s_i}))$ fall one level":}\\
 $x \leftarrow (\nu(T_{s_i}) - \mu'(T_{s_i}))$ \\
 $\mu'(s_i) \leftarrow \mu'(s_i) + x $ \% \textit{and simultaneously}\\
 $\mu'(r) \leftarrow \mu'(r) - x $.\\
 % Le coût de cette modification est $\ell_{T_i}(\nu(T_i) - \mu'(T_i)) = \ell_{T_i}|\mu'(T_i) -  \nu(T_i) |  $\\
$j \leftarrow \min\{k: \sum_{n=1}^{k}\pi_{r,n} \geq x \}$ \\
\textbf{for } $k < j$ 
\begin{myEnumerate}
\item $\pi_{s_i,k} \leftarrow \pi_{s_i,k} + \pi_{r,k}$
\end{myEnumerate}
\textbf{end for}\\
$\pi_{s_i,j} \leftarrow \pi_{s_i,j} + (x - \sum_{n=1}^{j-1} \pi_{r,n})$\\
$\pi_{r,j} \leftarrow \pi_{s_i,j} - (x - \sum_{n=1}^{j-1} \pi_{r,n})$\\
\textbf{for} $k < j$
\begin{myEnumerate}
\item $\pi_{r,k} \leftarrow 0$
\end{myEnumerate}
\textbf{end for}

\end{myEnumerate}
\textbf{end if}\\
$M \leftarrow M+1$
\end{myEnumerate}
\textbf{end for}
\end{myEnumerate}
\textbf{end if}
\end{myEnumerate}
\textbf{end for}
\end{myEnumerate}
\textbf{end for}

We must now prove that the algorithm works as intended, that is $\pi'$ is always a coupling between $\mu$ and $\mu'$, and  when the algorithm terminates we have $\mu' = \nu$ and the cost of $\pi'$ is $\sum_{e \in E} w_e \mid \mu(T_e) - \nu(T_e) \mid$.\\
 
\bprr
A probability measure on a tree $T$ is determined by the measure attributed to all subtrees $T_e$. To see that $\mu' = \nu$ when the algorithm terminates, we thus show that $\mu'(T_e) = \nu(T_e)$ for all subtree $T_e$:
\begin{itemize}
\item If $\mu(T_e) = \nu(T_e)$ then neither phase (1) nor (2) modifies $\mu'(T_e) = \mu(T_e) = \nu(T_e)$ (even though the distribution may vary).

\item If $\mu(T_e) > \nu(T_e)$ then phase (1) removes the adequate quantity of mass from $\mu'(T_e)$ so that once phase (1) is over we have $\mu'(T_e) = \nu(T_e)$. Then phase (2) does not change the quantity $\mu'(T_e) = \nu(T_e)$ (even though it could change the distribution on that subtree). 

\item If $\mu(T_e) < \nu(T_e)$ phase (1) does not change the quantity $\mu'(T_e) = \mu(T_e)$ (even though it could change the distribution on that subtree). We write $\mu^N$ for the measure $\mu'$ after all subtrees whose root is at depth $N$ have been treated by phase (2) ($N$ going from $0$ to the deepest level in the tree $-1$). Then we proceed by induction on $N$, assuming $e^+$ is at depth $N$. The initial step consists in seeing that $\mu^{N=0}$, the measure $\mu'$ just after phase (1), is a probability measure on $T$; $\nu$ being one too it follows $\mu^{N=0}(T) = \nu(T) = 1$. For the induction step, we write $e^- = v_1,...,v_m$ the children of $e^+$ and assume that (induction hypothesis) for all $i = 1,...,m$:
\begin{align*}
\mu^N(T_{e^+}) = \nu(T_{e^+}) &= \mu^N(e^+) + \sum_{i=1}^m \mu^N(T_{v_i})\\
 &=\nu(e^+) + \sum_{i=1}^m \nu(T_{v_i}).
\end{align*}
Since phase (1) is over $\mu^N(T_{v_i}) \leq \nu(T_{v_i})$, hence $\mu^N(e^+) \geq \nu(e^+)$. Then, phase (2) of the algorithm modifies:
$$\mu^{N+1}(v_i) = \mu^{N}(v_i) + (\nu(T_{v_i}) - \mu^N(T_{v_i})).$$
And:
\begin{align*}
\mu^{N+1}(T_{v_i}) &= \mu^{N}(T_{v_i}) - \mu^{N}(v_i) + \mu^{N+1}(v_i)\\
 &=  \mu^{N}(T_{v_i}) - \mu^{N}(v_i) + \mu^{N}(v_i) + (\nu(T_{v_i}) - \mu^N(T_{v_i}))\\
 &= \nu(T_{v_i}).
\end{align*}
Then we have the desired fact for $i = 1$.
\end{itemize}
Eventually when the algorithm stops $\mu' = \nu$.

$\mu'$ and $\pi'$ are modified only during loops (*) and (**), we write $\mu^M$ and $\pi^M$ for the values of $\mu'$ and $\pi'$ after $M$ rounds through loops (*) or (**). Then $ \pi^M = (\pi^M (x,y))_{x,y \in V} := (\pi^M_{x,y})_{x,y \in V} $ is a coupling between $\mu$ and $\mu^M$: just after initialization it is clear that $\pi' = \pi^0$ is a coupling between $\mu$ and $\mu^0 = \mu$, if follows by induction that $\pi^M$ is a coupling between $\mu$ and $\mu^M$ (treating separately the case where moving from $M$ to $M+1$ is done during phase (1) and the case where this move is done during phase (2)). About the cost of the coupling, if moving from $M$ to $M+1$ is done during phase (1), in loop (*) we have:
\begin{align*}
\sum_{x,i} d(x,i) \pi^{M+1}_{x,i} &= \sum_x \sum_i d(x,i) \pi^{M+1}_{x,i}\\
 &= \sum_i d(s,i) \pi^{M+1}_{s,i} + \sum_i d(r,i)  \pi^{M+1}_{r,i} + \sum_{x \neq r,s} \sum_i d(x,i) \pi^M_{x,i}\\
 &= \sum_{i< j} d(s,i)(\pi_{s,i}^M + \pi_{r,i}^M) + d(s,j) \bigg(\pi_{s,j}^M + \bigg(x - \sum_{n=1}^{j-1} \pi_{r,n}^M \bigg)\bigg) + \sum_{i > j} d(s,i) \pi_{s,i}^M\\
 & \ + \sum_{i< j} d(r,i)\cdot 0 + d(r,j) \bigg(\pi_{r,j}^M - \bigg(x - \sum_{n=1}^{j-1} \pi_{r,n}^M \bigg)\bigg) + \sum_{i > j} d(r,i) \pi_{s,i}^M\\
 & \ + \sum_{x \neq r,s} \sum_i d(x,i) \pi^M_{x,i}\\
 &= \sum_{x,i}d(x,i) \pi_{x,i}^M + (d(s,j) - d(r,j))\bigg(x - \sum_{n=1}^{j-1} \pi_{r,n}^M \bigg) + \sum_{i <j} (d(s,i) - d(r,i))\pi_{r,i}^M\\
 &= \sum_{x,i}d(x,i) \pi_{x,i}^M +  x \cdot (d(s,j) - d(r,j)) \\
 &= \sum_{x,i}d(x,i) \pi_{x,i}^M + x \cdot d(s,r) .
\end{align*}
Where we used that $d(s,j) - d(r,j) = d(s,r)$ and that every $i$ such that $\pi_{r,i}^M$ contributes to the sum $\sum_{i < j}(d(s,j)-d(r,j))\pi_{r,i}^M$ is such that $d(s,i) - d(r,i) = d(s,j)-d(r,j) = d(s,r) \geq 0$: Each vertex $i$ such that $\pi_{r,i}^M \neq 0$ is (non strictly) below $r$ in the rooted-tree  (since $\pi_{r,i}^M$ is the mass $\mu^M$ in $r$ coming from $i$; we let the reader check it formally). Then those vertices $i$ such that $\pi_{r,i}^M$ contributes to the sum $\sum_{i< j} (d(s,i)-d(r,i))\pi_{r,i}^M$ are  (non strictly) below  $r$ and for those we have $d(s,i) \geq d(r,i)$ and then $d(s,i) - d(r,i) = d(s,r)$ since $s$ is the father of $r$. By definition of $j$, $\pi_{r,j}^M \neq 0$ and thus $j$ is (non strictly) below $r$, hence $d(s,j) - d(r,j) = d(s,r) \geq 0$. If moving from $M$ to $M+1$ is done during phase (2), in loop (**), we conclude similarly that the cost of the coupling is increased by $x \cdot d(s_i,r)$. During phase (1) excess measure is always brought up one level at the time, in the loop (*) we thus always have $x = \mu'(T) - \nu(T) =  \mu(T) - \nu(T)$. And phase (1) brings up excess measure exactly from those subtrees $T_e$ with $\mu(T_e) > \nu(T_e)$. During phase (2) measure is always brought down one level at the time, in the loop (**) we thus always have $ x = \mu'(T_i) - \nu(T_i) = \mu(T_i) - \nu(T_i)$. And phase (2) brings down adequate quantity of measure exactly in those subtrees $T_e$ with $\mu(T_e) < \nu(T_e)$. Since just after initialization $\pi'$ has null cost, the cost of it at the end of the algorithm is thus $\sum_{e \in E} w_e | \mu(T_e) - \nu(T_e) |$.
\epr

\begin{Rem} Let $(X,d)$ be a Polish metric space. For $\mu,\nu\in P_1(X)$, we have from the {\it Kantorovich-Rubinstein duality}:
$$Wa(\mu,\nu)=\sup\Big\{\int_X f(x)\,d\mu(x)-\int_X f(x)\,d\nu(x)\Big\},$$
where the supremum is taken over all 1-Lipschitz functions $f$: see Theorem 1.3 in \cite{Villa}; see also \cite{Edwards} for a short proof. We observe that, for a finite metric tree, our second proof of Theorem \ref{WassTrees} does not appeal to Kantorovich-Rubinstein duality (in contrast e.g. with the proof in \cite{EvMat}).
\end{Rem}

\section{Appendix: $\sigma$-algebras on real trees}

Let $(T,d)$ be a real tree. Apart from Godard's construction from \cite{God} of the $\sigma$-algebra $\mathcal{G}$ recalled above, we are aware of other constructions of $\sigma$-algebras on $T$ and of corresponding length measures: 
\begin{itemize}
\item The $\sigma$-algebra $\mathcal{S}$ generated by segments, see \cite{Valette}.
\item The Borel $\sigma$-algebra $\mathcal{B}$ generated by open subsets, see \cite{EPW} for compact real trees, then for locally compact real trees in \cite{AEW}.
\end{itemize} 
All these constructions have in common that the length measure of a segment $[x,y]$ is exactly $d(x,y)$. In order to clarify the relation between $\mathcal{S}, \mathcal{B}$ and $\mathcal{G}$, we also introduce the $\sigma$-algebras $\mathcal{B}_0$ generated by open balls (so that $\mathcal{B}_0\subset\mathcal{B}$) and $\overline{\mathcal{S}}$ obtained by completing $\mathcal{S}$ with respect to $\lambda$-negligible subsets.

The following proposition explains our choice to work with Godard's $\sigma$-algebra $\mathcal{G}$.

\begin{Prop} Let $T$ be a real tree. 
\begin{enumerate}
\item We have $\mathcal{S}\subset\mathcal{B}_0\subset\mathcal{G}$ and $\overline{\mathcal{S}}\subset\mathcal{G}$.
\item If $T$ is separable then $\mathcal{B}_0=\mathcal{B}$ and $\overline{\mathcal{S}}=\mathcal{G}$
\end{enumerate} 
\end{Prop}

\bprr \begin{enumerate}
\item To show that $\mathcal{S}\subset\mathcal{B}_0$, fix $x,y\in T$ and let $(z_n)_{n>0}$ be a dense sequence in $[x,y]$. Then the equality
$$[x,y]=\bigcap_{m\geq 1}\Big(\bigcup_{n>0} B(z_n, 1/m) \Big)$$
shows that $[x,y]\in\mathcal{B}_0$. Now, let $B$ be an open ball in $T$. For any $x,y\in T$, the intersection $B\cap [x,y]$ is convex in $[x,y]$, so it is a sub-interval in $[x,y]$. In particular $B\cap [x,y]$ is Lebesgue-measurable in $[x,y]$, so $B\in\mathcal{G}$.

The inclusion $\overline{\mathcal{S}}\subset\mathcal{G}$ follows from the fact that $\mathcal{G}$ is complete, as can be seen from the definitions.

\item The equality $\mathcal{B}_0=\mathcal{B}$ holds in every separable metric space (any open set being then a countable union of open balls).

To prove the inclusion $\mathcal{G}\subset\overline{\mathcal{S}}$, we consider the subset $T^0=:\cup_{x,y\in T}]x,y[$ and its complement $L=T\setminus T^0$: the latter is the set of {\it leaves} of $T$. For every segment $[x,y]$ we have $L\cap [x,y]\subset \{x,y\}$, so that $L\in \mathcal{G}$; moreover $\lambda(L)=0$ by equation (\ref{Godmeasure}).

%Fix now $A\in\mathcal{G}$. Then $A\cap L$ is $\lambda$-negligible, so $A\cap L\in\overline{\mathcal{S}}$. To show $A\cap T^0\in \overline{\mathcal{S}}$, we use separability of $T$: let $D$ be a countable dense subset, then it is easy to see that $T^0=\cup_{x,y\in D}]x,y[$, so $A\cap T^0$ is the countable union of the $A\cap ]x,y[$'s, with $x,y\in D$. But in a segment, the $\sigma$-algebra of Lebesgue-measurable subsets is the completion of the $\sigma$-algebra generated by sub-intervals, i.e. $A\cap ]x,y[\in \overline{\mathcal{S}}$. This concludes the proof.

Then we take $A \in \mathcal{G}$. To show $A \in \overline{\mathcal{S}}$ we use separability: let $D$ be a countable subset of $T$. It is easy to see that $T^\circ = \cup_{x,y \in D} ]x,y[$ which implies that $T^\circ$ is $\overline{\mathcal{S}}$-measurable, as well as its complement $L$. On the one hand $A \cap L \subset L$ and $\lambda(L) =0$, then $A \cap L$ is $\lambda$-negligible and thus $\overline{\mathcal{S}}$-measurable. On the other hand $A\cap T^\circ \in \overline{\mathcal{S}}$ since $A \cap ]x,y[ \in \overline{\mathcal{S}}$ for all $x,y \in T$ because the $\sigma$-algebra of Lebesgue-measurable subsets is the completion of the $\sigma$-algebra generated by sub-intervals i.e. $A \cap ]x,y[ \in \overline{\mathcal{S}}$. This concludes the proof.

\epr
\end{enumerate}

\noindent
Authors' address:\\
Institut de Math\'ematiques\\
Universit\'e de Neuch\^atel\\
Unimail\\
11 Rue Emile Argand\\
CH-2000 Neuch\^atel - SWITZERLAND

\medskip
\noindent
maxmatheyp@gmail.com\\
alain.valette@unine.ch


\begin{thebibliography}{CCJJV}

\bibitem[AEW13]{AEW} S. {\sc Athreya}, M. {\sc Eckhoff} and A. {\sc Winter},
\newblock {\em Brownian motion on $\R$-trees},
\newblock Transactions of the American Mathematical Society 365 (2013), 3115-3150.

\bibitem[BMSZ20]{BMSZ} F. {\sc Baudier}, P. {\sc Motakis}, G. {\sc Schlumprecht} and A. {\sc Zsak},
\newblock {\em Stochastic approximation of lamplighter metrics},
\newblock Preprint, arXiv:2003.06093

\bibitem[Bo85]{Bour} J. {\sc Bourgain}, 
\newblock {\em On Lipschitz embedding of finite metric spaces in Hilbert space},
\newblock Israel Journal of Mathematics, 52, 46-52, 1985.

\bibitem[CC$^+$98]{CCGGP} M. {\sc Charikar}, C. {\sc Chekuri}, A. {\sc Goel}, S. {\sc Guha}, and S. {\sc Plotkin},
\newblock {\em Approximating a finite metric by a small number of tree metrics},
\newblock Proc. Symposium Foundations of Computer Science, 1998.

\bibitem[Ch02]{Cha} M.S. {\sc Charikar},
\newblock {\em Similarity estimation techniques from rounding algorithms}
\newblock STOC '02: Proceedings of the thirty-fourth annual ACM symposium on Theory of computing, May 2002, 380-388.

\bibitem[DL97]{DeLa} M.M. {\sc Deza} and M. {\sc Laurent},
\newblock {\em Geometry of cuts and metrics},
\newblock Springer, Algorithms and combinatorics 15, 1997.

\bibitem[Ed10]{Edwards} D.A. {\sc Edwards},
\newblock {\em A simple proof in Monge-Kantorovich duality theory},
\newblock Studia Math. 200 (2010), 67-77.

\bibitem[EPW06]{EPW} S.N. {\sc Evans}, J. {\sc Pitman} and A. {\sc Winter},
\newblock {\em Rayleigh processes, real trees, and root growth with re-grafting},
\newblock Probab. Theory Relat. Fields 134 (2006), 81-126.

\bibitem[EM12]{EvMat} S.N. {\sc Evans} and F.A. {\sc Matsen}
\newblock {\em The phylogenetic Kantorovich-Rubinstein metric for environmental sequence samples},
\newblock J. R. Stat. Soc. Series B Stat. Methodol. 74(3) (2012), 569-592.

\bibitem[FRT04]{FaRaTa} J. {\sc Fakcharoenphol}, S. {\sc Rao} and K. {\sc Talwar},
\newblock {\em A tight bound on approximating arbitrary metrics by tree metrics},
\newblock Journal of Computer and System Sciences 69 (2004) 485-497.


\bibitem[Go10]{God} A. {\sc Godard}.
\newblock {\em Tree metrics and their Lipschitz-free spaces}
\newblock Proc. Amer. Math. Soc. 138 (2010), 4311-4320.

\bibitem[Ha79]{Haagerup}   U. {\sc Haagerup},
\newblock {\em An example of a nonnuclear $C^*$-algebra, which has the metric approximation property},
\newblock Invent. Math. 50 (1978/79), 279-293.

\bibitem[IT03]{InTha} P. {\sc Indyk} and N. {\sc Thaper}.
\newblock {\em Fast image retrieval via embeddings}. 
\newblock In 3rd International Workshop on Statistical and Computational Theories of Vision (at ICCV), 2003.

\bibitem[KT02]{KleTar} J. {\sc Kleinberg} and E. {\sc Tardos},
\newblock {\em Approximation Algorithms for Classification Problems with Pairwise Relationships: Metric Labeling and Markov Random Fields},
\newblock Journal of the ACM 49 (2002), 616-639.

\bibitem[Kl15]{Kloeck} B.R. {\sc Kloeckner},
\newblock {\em A geometric study of Wasserstein spaces: ultrametrics},
\newblock Mathematika 61 (2015), 162-178. 

\bibitem[Na18]{NaorICM} A. {\sc Naor}
\newblock {\em Metric dimension reduction: a snapshot of the Ribe program}. 
\newblock Proceedings of the ICM-Rio de Janeiro 2018. Vol. I. Plenary lectures, 759-837, World Sci. Publ., Hackensack, NJ, 2018. 

\bibitem[NS07]{NaSche} A. {\sc Naor} and G. {\sc Schechtman},
\newblock {\em Planar earthmover is not in $L^1$}. 
\newblock SIAM J. Comput. 37 (2007), 804-826. 


\bibitem[OO19]{OsOs} S. {\sc Ostrovska} and M.I. {\sc Ostrovskii},
\newblock {\em Generalized transportation cost spaces}
\newblock, Preprint, arXiv:1902.10334

\bibitem[Os13]{Ostro} M.I. {\sc Ostrovskii},
\newblock {\em Metric embeddings - bilipschitz and coarse embeddings into Banach spaces},
\newblock De Gruyter, Studies in Math. 49, 2013.

\bibitem[Sa15]{Santa} F {\sc Santambrogio},
\newblock {\em Optimal Transport for Applied Mathematicians},
\newblock Birkh\"auser, Progress in Nonlinear Differential Equations and Their Applications, 2015.

\bibitem[Va90]{Valette} A. {\sc Valette},
\newblock {\em Les repr\'esentations uniform\'ement born\'ees associ\'ees à un arbre r\'eel}, 
\newblock Bull. Soc. Math. Belg. S\'er. A 42 (1990), 747-760. 

\bibitem[Vi15]{Villa} C. {\sc Villani},
\newblock {\em Topics in Optimal Transportation},
\newblock AMS, Graduate Studies in Mathematics Vol. 58, 2003.

\bibitem[We99]{Weaver} N. {\sc Weaver},
\newblock {\em Lipschitz algebras}. 
\newblock World Scientific Publishing Co., Inc., River Edge, NJ, 1999. xiv+223 pp.

\end{thebibliography}
\end{document}